\documentclass{amsart}
\usepackage{amssymb,graphics}

\theoremstyle{plain}
\newtheorem{theorem}{Theorem}

\newtheorem{proposition}[theorem]{Proposition}

\theoremstyle{definition}
\newtheorem{definition}{Definition}
\newtheorem{example}{Example}

\newcommand\Z{{\mathbb Z}}
\newcommand\Aut{{\mathsf{Aut}}}
\newcommand\Sym{{\mathsf{Sym}}}

\newcommand\St{{\mathsf{St}}}
\newcommand\GL{{\mathsf{GL}}}
\newcommand\tree{{\mathcal T}}
\newcommand\btree{{\partial \mathcal T}}
\newcommand\A{{(Q,X,\rho,\tau)}}

\title{Some solvable automaton groups}

\author{Laurent Bartholdi}
\address{\'Ecole Polytechnique F\'ed\'erale de Lausanne (EPFL), Institut de math\'ematiques B (IMB), CH-1015 Lausanne, Switzerland}
\email{laurent.bartholdi@epfl.ch}
\author{Zoran \v Suni\'k}
\address{Department of Mathematics, Texas A\&M University, College Station, TX 77843-3368, USA}
\email{sunik@math.tamu.edu}

\date{March 14, 2006}

\begin{document}
\begin{abstract}
It is shown that certain ascending HNN extensions of free abelian
groups of finite rank, as well as various lamplighter groups, can
be realized as automaton groups, i.e., can be given a self-similar
structure. This includes the solvable Baumslag-Solitar groups
$BS(1,m)$, for $m \neq \pm 1$.

In addition, it is shown that, for any relatively prime integers $m,
n \geq 2$, the pair of Baumslag-Solitar groups $BS(1,m)$ and
$BS(1,n)$ can be realized by a pair of dual automata. The examples
are then used to illustrate more general connections between Schreier
graphs, composition of automata and dual automata.
\end{abstract}
\maketitle


Groups generated by automata appeared already in the 1950's. Among
the pioneering works we mention Horejs~\cite{horejs:automata} and
Aleshin~\cite{aleshin:burnside}. Important examples appeared
later, in particular the well known examples of infinite
residually finite torsion groups, and groups of intermediate
growth constructed by Grigorchuk
in~\cite{grigorchuk:burnside,grigorchuk:growth}. Many groups were
then shown to belong to that class; in particular linear groups
over $\Z$~\cite{brunner-s:glnz}.

The set of all transformations generated by finite automata over a
fixed finite alphabet form a group, denoted $\mathcal F$. It is
not known which solvable groups appear as subgroups of $\mathcal
F$, i.e., appear as groups generated by finite automata. Progress
in this direction has been achieved in the works of Sidki and
Brunner~\cite{brunner-s:wreath,sidki:adding,sidki:wreathing}.

In this note, we are interested in (solvable) groups that are
generated by all the states of a single finite automaton. Such
groups are called automaton groups. The special interest in this
more restricted setting is justified by the self-similarity
structure that is apparent as soon as a group is realized as an
automaton group.

The purpose of this note is twofold. We go over some well known
notions and constructions (automaton groups, inversion,
composition) as well as some less known (dual automata). At the
same time, we realize some solvable groups as automaton groups
(thus giving them self-similar structure) and use them to
illustrate the introduced notions.

For example, we show that, for any $n$ coprime to $m$, the
solvable Baumslag-Solitar groups
\[ BS(1,m) = \langle \ a,t \ |\ tat^{-1} = a^m \ \rangle \]
belong to the class of automaton groups on a $n$-letter alphabet.
The automata that describe them are related to multiplication by
$m$ and addition in base $n$.

Similar constructions, corresponding to multiplication by linear
polynomials over the finite ring $\Z/nZ$, lead to ``lamplighter
groups'', i.e.\ the groups
\[ L_n = (\Z/n\Z) \wr \Z = \langle \ a,t \ |\ a^n=[a,t^iat^{-i}]=1\
  \forall i\in\Z\ \rangle. \]

The above considerations are then extended to the
multi-dimensional case. Namely, for any $d \geq 1$ and any $d
\times d$ matrix $M$ of infinite order and determinant $m$
relatively prime to $n$, the ascending HNN extension of the free
abelian group of rank $d$ by the endomorphism defined by $M$ can
also be realized by a finite automaton. In this case, the
automaton corresponds to multiplication by the matrix $M$ in the
free $d$-dimensional module over $n$-adic integers.

Similarly, automata corresponding to multiplication by monic
invertible polynomials of degree $d$ over the finite ring $Z/n\Z$
lead to construction of lamplighter groups of the form
$L_{n,d}=(\Z/nZ)^d \wr \Z$.

The lamplighter group $L_2$ was realized by a 2-state automaton by
Grigorchuk and {\. Z}uk in~\cite{grigorchuk-z:l2}. During the
preparation of this manuscript the authors have learned that Silva
and Steinberg have also constructed various lamplighter groups by
using finite automata in~\cite{silva-s:lamplighter}. Their
construction is based on the so called reset automata, for which
the alphabet and the set of states are usually the same. Thus the
realization of $L_{n,d}$ can be done by an $n^d$-state
$n^d$-letter reset automaton. Our results show that the
$n^d$-state automaton $A_{1+t^d}$ acting on the $n$-ary rotted
tree also defines $L_{n,d}$. However, Silva and Steinberg point
out that the construction involving reset automata is essentially
the simplest in terms of Krohn-Rhodes theory.


\section*{Tree automorphisms}
Let $X$ be a finite alphabet. The set of finite words $X^*$ over $X$
has a structure of a rooted labelled $n$-ary tree, denoted $\tree(X)$
or sometimes simply $\tree$. The empty word $\emptyset$ is the root
of the tree and the words of length $k$ constitute the vertices on
the level $k$, denoted $L_k$, in the tree. A vertex $u$ on  level $k$
is a neighbor to a vertex $v$ on level $k+1$ if and only if $v=ux$
for some letter $x \in X$. A word $u$ is a prefix of a word $v$ if
and only if there exists a word $w$ such that $v=uw$. This is
equivalent to the condition that $u$ is a vertex on the unique path
from the root to $v$. The group of automorphisms $\Aut(\tree)$ of the
tree $\tree$ consists of all permutations of $X^*$ that preserve the
structure of the tree. Such permutations must preserve the root,
since the root is the only vertex of degree $n$, must preserve the
levels, since the distance to the root must be preserved, and must
preserve the prefix relation, since paths are mapped to paths. The
group $\Aut(\tree)$  consists precisely of those permutations of
$X^*$ that preserve the prefix relation. The boundary $\btree$ of
$\tree$ is a metric space $(X^\omega,d)$ whose elements are the
infinite rays in $\tree$ starting at the root (right infinite words
over $X$). The distance $d$ between two distinct rays $r$ and $\ell$
in $\btree$ is defined by $d(r,\ell) = 2^{-|r \wedge \ell|}$, where
$|r \wedge \ell|$ denotes the length of the longest common prefix $r
\wedge \ell$ of $r$ and $\ell$. There is a canonical isomorphism
between $\Aut(\tree)$ and the group of isometries of $\btree$. Given
an isometry $\overline{f}$ of $\btree$ define an automorphism $f$ of
$\tree$ as follows. For a word $w$ of length $k$ define $f(w)$ to be
the prefix of length $k$ of the image $\overline{f}(r)$ of any ray
$r$ that has $w$ as a prefix. We find it useful to sometimes switch
back and forth between these two interpretations of $\Aut(\tree)$,
i.e., we may define tree automorphisms by defining the action on
infinite words.

The $|X|$ trees hanging below the root are canonically
isomorphic to $\tree$. Thus the stabilizer $\St(L_1)$ of the
first level in $\tree$ is canonically isomorphic to
$\Aut(\tree)^X$. The symmetric group $\Sym(X)$ on $X$ canonically
embeds in $\Aut(\tree)$ as the group of rooted tree automorphisms
defined by
\[ \rho(xw) = \rho(x)w, \]
for $\rho$ in $\Sym(X)$, $x$ a letter in $X$ and $w$ a word over
$X$. The stabilizer $\St(L_1)=\Aut(\tree)^X$ is normal in
$\Aut(\tree)$ and the  group of rooted tree automorphisms is its
transversal, leading to the permutational wreath product
decomposition
\[ \Aut(\tree) = \Aut(\tree)^X \rtimes \Sym(X)  = \Aut(\tree) \wr \Sym(X).\]
The symmetric group $\Sym(X)$ acts on the right of $\Aut(\tree)^X$
by
\[ (f^\rho)_x = f_{\rho(x)} \]
for $\rho \in \Sym(X)$ and $f \in \Aut(\tree)^X$ (here $f_x$
denotes the automorphism in $\Aut(\tree)$ that is at the
$x$-component of $f$). Each tree automorphism $f$ can be written
uniquely as
\[ f = \rho_f (f_x)_{x \in X} \]
where $f_{x}$, called the section of $f$ at $x$, is a
tree automorphism corresponding to the way $f$ acts on the
subtree $\tree_x$ consisting of the words that start in $x$, and
$\rho_f$, called the root permutation of $f$, is a
permutation of $X$ corresponding to the way $f$ permutes the
$|X|$ subtrees below the root. The root permutation $\rho_f$ of
$X$ and the sections automorphisms $f_x$, $x \in X$, are
determined uniquely from the equalities
\begin{equation} \label{fxw}
f(xw) = \rho_f(x)f_x(w),
\end{equation}
for $x$ a letter in $X$ and $w$ a word over $X$. Since $\rho_f$
is just the restriction of $f$ on $X$ we may write
\begin{equation} \label{fxw2}
f(xw) = f(x)f_x(w),
\end{equation}
for $x$ a letter in $X$ and $w$ a word over $X$. The composition
of two tree automorphisms $f$ and $g$ is an automorphism,
denoted $fg$, with
\begin{equation}\label{composition}
 \rho_{fg}=\rho_f \rho_g \qquad\text{and}\qquad
 (fg)_x = f_{g(x)}g_x,
\end{equation}
for $x \in X$. For the inverse $f^{-1}$ we have
\begin{equation}\label{inverse}
 \rho_{f^{-1}} = \rho_f^{-1} \qquad\text{and}\qquad
 (f^{-1})_x = (f_{f^{-1}(x)})^{-1},
\end{equation}
for $x \in X$.


\section*{Automata as tree automorphisms}
We now define special kind of tree automorphisms, defined by
finite automata. A good reference for these constructions
is~\cite{grigorchuk-n-s:automata}.

\begin{definition}
A \emph{finite synchronous transducer} is a quadruple
\[ A = \A, \]
where $Q$ is a finite set whose elements are called \emph{states},
$X$ is a finite set called the \emph{alphabet} of $A$ and whose
elements are called \emph{letters}, and the functions
\[ \rho: Q \times X \to X \qquad\text{and}\qquad
   \tau: Q \times X \to Q \]
are called the \emph{rewriting} and the \emph{transition}
functions of $A$.
\end{definition}

We refer to finite synchronous transducers simply by calling them
automata. The rewriting and transition function define a recursive
way in which every state of the automaton $A=\A$ rewrites the words
over $X$. When the automaton is in state $q$ and is faced with the
input word $xw$ it rewrites the input letter $x$ into the output
letter $\rho(q,x)$ and changes its state to $\tau(q,x)$, which state
then handles $w$, i.e., the rest of the input. In other words, the
domains of the rewriting and transition functions are extended (in
the second variable) to arbitrary words by
\begin{gather}
 \rho(q,xw) =  \rho(q,x)\rho(\tau(q,x),w), \label{rhoqxw}\\
 \tau(q,xw) =  \tau(\tau(q,x),w) \label{tauqxw}.
\end{gather}

\begin{definition}
An automaton $A=\A$ is invertible if, for each state $q$ in $Q$,
the restriction $\rho_q: X \to X$, defined by
$\rho_q(x)=\rho(q,x)$, is a permutation.
\end{definition}

Consider an invertible automaton $A=\A$. By introducing notation
$\rho(q,w)=q(w)$ and $\tau(q,w)=q_w$, the equalities (\ref{rhoqxw})
and (\ref{tauqxw}) can be rewritten (compare to (\ref{fxw}) and
(\ref{fxw2})) as
\begin{gather*}
 q(xw) =  \rho_q(x)q_x(w) = q(x)q_x(w), \\
 q_{xw} =  (q_x)_w.
\end{gather*}

Each state $q$ of an invertible automaton defines an automorphism,
also denoted $q$, of the regular rooted $|X|$-ary tree. Note that the
notation $\rho_q$ and $q_x$ is consistent with the earlier notation
used for tree automorphisms, since $\rho_q$ is indeed the root
permutation of $X$ induced by the automorphism $q$ and $q_x$ is the
section of $q$ at $x$.

\begin{example}\label{f}
Let $X$ be a finite set and $f:X^{d+1} \to X$ an arbitrary map.
Define an automaton $A_f=(X^d,X,\rho,\tau)$, where $\rho:X^d
\times X \to X$ and $\tau:X^d \times X \to X^d$ are given by
\[ \rho((x_1,\dots,x_d),x) = f(x_1,\dots,x_d,x) \quad\text{and}\quad
   \tau((x_1,\dots,x_d),x) = (x_2,\dots,x_d,x), \]
respectively. It follows directly from the definition that if,
for all $d$-tuples ${\bf y} \in X^d$, the restriction $f_{\bf
y}: X \to X$ given by $x \mapsto f({\bf y},x)$ is a permutation,
the automaton $A_f$ is invertible. The tree automorphism defined
by the state ${\bf y}=(y_1,\dots,y_d) \in X^d$ is given by
\[ {\bf y}(x_1x_2x_3\dots) = f(y_1,\dots,y_d,x_1)f(y_2,\dots,y_d,x_1,x_2)
   f(y_3,\dots,y_d,x_1,x_2,x_3)\dots .\]
Note that only the first $d$ symbols of the output depend on the
state ${\bf y}$.

As a more special example, let $X$ be the finite ring $X=R=\Z/n\Z$ and
let $g=a_0 + a_1t+ \dots + a_dt^d$ be a monic polynomial of degree $d
\geq 1$, which is invertible in the power series ring $R[[t]]$ (thus
we assume that $a_0$ is invertible in $R$ and $a_d=1$). Consider the
function $f:X^{d+1} \to X$ given by $f(x_0,x_1,\dots,x_d) =
a_dx_0+a_{d-1}x_1+ \dots + a_0x_d$. Then the automaton $A_f$, which
we also denote by $A_g$, is invertible. In particular, when $g=1+t$,
the rewriting and the transition functions of $A_{1+t}$ are given by
\[ \rho(y,x) = y + x \qquad\text{and}\qquad\tau(y,x) = x.\]
\end{example}

\begin{example}
For $a$ an integer and $b$ a positive integer, denote by $a \boxdot
b$ and $a \div b$ the remainder and the quotient obtained when $a$ is
divided by $b$.

For positive and relatively prime integers $m$ and $n$ define
the automaton
\[ S_{m,n} = (S,X,\rho,\tau) \]
where $S=\{s_0,\dots,s_{m-1}\}$, $X = \{x_0,\dots,x_{n-1}\}$,
and $\rho : S \times X \to X$ and $\tau: S \times X \to S$ are
given by
\[ \rho(s_i,x_j)= x_{(mj+i) \boxdot n} \qquad\text{and}\qquad
   \tau(s_i,x_j) = s_{(mj+i) \div n}, \]
respectively. The automaton $S_{m,n}$ is invertible. This is because
$m$ is invertible in $\Z/n\Z$ and therefore the map $j \mapsto mj+i$
is a permutation of $\Z/n\Z$.
\end{example}

We mention yet another way to think of tree automorphisms defined by
finite invertible automata. Let $Q$ be a finite set of symbols and
let $\rho_q$, for $q \in Q$, be a permutation of the alphabet
$X=\{x_1,\dots,x_n\}$. Consider the system of $|Q|$ equations
\[ q = \rho_q (q_1,\dots,q_n), \qquad\text{for } q \in Q \]
where $q_i \in Q$, for all $q$ and $i$. Such a system has a
unique solution in $\Aut(\tree)$ for all $q \in Q$, such that
$q_i \in Q$ is the section of $q$ at $x_i$ and $\rho_q$ is the
root permutation of $q$. The rewriting and the transition
functions in the automaton $A=\A$ that corresponds to the above
system of equations are given by
\[ \rho(q,x_i)=\rho_q(x_i) \qquad\text{and}\qquad \tau(q,x_i)=q_i, \]
for $q \in Q$ and $x \in X$.

\begin{figure}[!hbt]
\begin{center}
\includegraphics{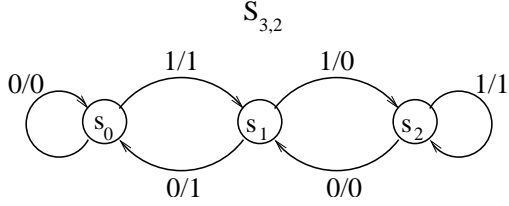}
\end{center}
\caption{The automaton $S_{3,2}$\label{s32}}
\end{figure}

An automaton $A=\A$ is usually depicted by a labelled directed graph
$\Gamma(A)$, where the set of vertices of $\Gamma(A)$ is $Q$ and a
directed edge from $q$ to $p$ labelled by $x|y$
\[ q \overset{x|q(x)} \longrightarrow q_x \]
exists in $\Gamma(A)$ if and only of $\rho(q,x)=y$ and $\tau(q,x)=p$,
i.e., $q(x)=y$ and $q_x=p$. Figure~\ref{s32} depicts the automaton
$S_3,2$ (with the agreement that $x_j=j$, for $j=0,1$.

Flipping every label $x|y$ to a label $y|x$ in the graph $\Gamma(A)$
of an invertible automaton leads to a graph of another invertible
automaton $\overline{A}$. Moreover, if $q$ is a vertex (state) in the
original graph (automaton) $\Gamma(A)$ then the corresponding vertex
(state) $\overline{q}$ in $\Gamma(\overline{A})$ defines the inverse
automorphism $q^{-1}$ of $q$ in $\Aut(\tree)$. Indeed, starting from
the state $q$ in $\Gamma(A)$ the automaton $A$ reads the word
$x_1x_2x_3\dots$ and outputs $y_1y_2y_3\dots$ while passing through
the states $q_{x_1}, q_{x_1x_2}, q_{x_1x_2x_3},\dots$. Starting from
the state $\overline{q}$, the automaton $\overline{A}$ reads the word
$y_1y_2y_3\dots$, follows the corresponding edges in
$\Gamma(\overline{A})$ and gives the output $x_1x_2x_3$ while passing
through the corresponding states $\overline{q_{x_1}},
\overline{q_{x_1x_2}}, \overline{q_{x_1x_2x_3}},\dots$. This simple
observation leads to the following definition.

\begin{definition}
For an invertible automaton $A=\A$, define the \emph{inverse}
automaton of $A$, denoted by $\overline A$, by
\[\overline A = (\overline Q,X, \overline \rho, \overline \tau) \]
where $\overline Q = \{\ \overline q \ |\  q \in Q \ \}$ is a
copy of the set $Q$, and $\overline \rho : \overline Q \times X
\to X$ and $\overline \tau: \overline Q \times X \to \overline
Q$ are given by
\[ \overline \rho(\overline q,x)=\rho_q^{-1}(x) \qquad\text{and}\qquad
   \overline \tau(\overline q,x) = \overline{\tau(q,\rho_q^{-1}(x))} . \]
\end{definition}

Note that the definition looks rather convoluted, even though
all we did is flip all the labels. Using the simplified
notation, we may write
\[ \overline \rho_{\overline q} = \rho_{q}^{-1} \quad\text{and}\quad
   \overline q_x = \overline{q_{q^{-1}(x)}}, \]
for a state $\overline q$ in $\overline Q$ and a letter $x$ in $X$,
which is compatible with the equalities~(\ref{inverse}).

\begin{example}\label{barA}
The automaton $\overline{A_{1+t}} =(\overline X,X,\overline \rho,
\overline \tau)$, where
\[ \overline \rho(\overline y,x) = (-y+x) \boxdot n \qquad\text{and}\qquad
   \overline \tau(\overline y,x) = \overline{(-y+x) \boxdot n}, \]
is the inverse of the automaton $A_{1+t}$. Figure~\ref{l2} depicts
the automaton $A_{1+t}$ and its inverse $\overline{A_{1+t}}$ in the
binary case (when $n=2$).

\begin{figure}[!hbt]
\begin{center}
\includegraphics{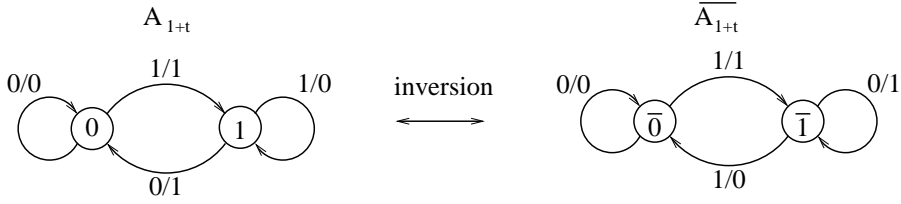}
\end{center}
\caption{The automaton $A_{1+t}$ and its inverse
$\overline{A_{1+t}}$\label{l2}}
\end{figure}
\end{example}

\begin{example}
The inverse of the automaton $S_{m,n}$ is the automaton
\[ \overline{S_{m,n}} = (\overline S,X, \overline \rho,\overline \tau) \]
where $\overline S=\{\overline{s_0},\dots,\overline{s_{m-1}}\}$, $X =
\{x_0,\dots,x_{n-1}\}$,and $\overline \rho: \overline S \times X \to
X$ and $\overline \tau : \overline S \times X \to \overline S$ are
given by
\[   \overline \rho(\overline s_i,x_j)= x_{(m'(j-i)) \boxdot n}
     \qquad\text{and}\qquad
     \overline \tau(\overline s_i,x_j) =
     \overline{s_{(m [m'(j-i) \boxdot n]+i) \div n}} , \]
respectively, and $m'$ is the multiplicative inverse of $m$ modulo
$n$. Indeed, if we denote the restriction $\rho_{s_i}$ by $\rho_i$,
then $\rho_i^{-1}$ is given by $x_j \mapsto x_{ m'(j-i) \boxdot n}$
and therefore
\[ \overline \rho(\overline s_i,x_j)= \rho_i^{-1}(x_j) = x_{m'(j-i) \boxdot n} \]
and
\[ \overline \tau(\overline s_i,x_j) =
   \overline{\tau(s_i,\rho_i^{-1}(x_j))} =
   \overline{\tau(s_i,x_{m'(j-i) \boxdot n})} =
   \overline{s_{(m [m'(j-i) \boxdot n] +i) \div n}}. \]
\end{example}

Occasionally we will need the notion of isomorphic automata.

\begin{definition}
Two automata $A_1=(Q_1,X_1,\rho_1,\tau_1)$ and
$A_2=(Q_2,X_2,\rho_2,\tau_2)$ are isomorphic if there exists a pair
of bijections $\alpha:Q_1 \to Q_2$ and $\beta:X_1 \to X_2$ that are
compatible with the transition and rewriting functions, i.e.,
\[ \alpha(\tau_1(q,x)) = \tau_2(\alpha(q),\beta(x)) \qquad\text{and}\qquad
   \beta(\rho_1(q,x)) = \rho_2(\alpha(q),\beta(x)), \]
for $q$ a state in $Q_1$ and $x$ a letter in $X_1$.
\end{definition}

Quite often an easy way to check if a pair of bijections is an
isomorphism between automata is to check if it is an isomorphism of
the corresponding labelled graphs representing the automata. In other
words, if $\alpha: Q_1 \to Q_2$ and $\beta: X_1 \to X_2$ are
bijections it suffices to check if, for every edge of the form
\[ q \overset{x|y}\longrightarrow p \]
in the graphical representation of $A_1$, there exists an edge of the
form
\[ \alpha(q) \overset{\beta(x)|\beta(y)}\longrightarrow \alpha(p) \]
in the graphical representation of $A_2$.

If the alphabet is fixed under an isomorphism, i.e., $A_1$ and $A_2$
share the same alphabet and $\beta$ is the identity map, the states
of $A_1$ and $A_2$ define the same set of automorphisms of the tree
$\tree(X)$. We write $A_1 \cong A_2$ for isomorphic automata. In case
the automorphism is canonical in some way we may write $A_1 = A_2$.


\section*{Automaton groups}
\begin{definition}
The group $G(A) = \langle \ \{\ q \ | \ q \in Q \ \} \ \rangle
\leq \Aut(\tree) $ generated by the states of an invertible
automaton $\A$ is called the \emph{group of the automaton} $A$.
Any group of automorphisms $G \leq \Aut(\tree)$ for which there
exists an automaton $A$ such that $G=G(A)$ is called an
\emph{automaton group}.
\end{definition}

Isomorphic automata generate isomorphic groups of tree automorphisms.
In case the alphabet is fixed under the automata isomorphism, the two
automaton groups are the same.

We reconsider now the automata from Example~\ref{f}.

\begin{proposition}
The group of the automaton $A_{1+t}$ is the lamplighter group
\[ L_n = (\oplus_\Z \Z/n\Z)\rtimes\Z  = (\Z/n\Z) \wr \Z, \]
where the action of $\Z$ on itself is by translations.
\end{proposition}
\begin{proof}
The infinite sequences over $X=\{0,1,\dots,n-1\}$ can be
interpreted as the elements of the power series ring
$R[[t]]$, where $R$ is the ring $\Z/n\Z$. Consider the functions
$\alpha,\mu:R[[t]] \to R[[t]]$ given by
\[ \alpha(p)=p+1 \qquad\text{and}\qquad \mu(p)=(1+t)p, \]
respectively. They both define automorphisms of the $n$-ary tree
$\tree(X)$. Let $G=\langle \alpha,\mu \rangle$. For $k \in \Z$,
\[ \mu^k\alpha\mu^{-k}(p) = \mu^k\alpha((1+t)^{-k}p) =
   \mu^k((1+t)^{-k}p +1) = p + (1+t)^k. \]
The automorphisms $\mu^k\alpha\mu^{-k}$, for $k \in \Z$,
have order $n$, commute, and generate the normal
closure $N$ of $\alpha$ in $G$, isomorphic to $\oplus_Z \Z/n\Z$.
On the other hand, the automorphism $\mu$ has infinite
order, which then shows that $N \cap \langle \mu \rangle = 1$.
Thus $G = (\oplus_Z \Z/n)\rtimes\Z$. Since conjugation by $\mu$
shifts the components in $N=\oplus_Z \Z/n\Z$ it is clear that $G
\cong L_n$.

It remains to be shown that the states of the automaton $A_{1+t}$
generate $G$. In order to avoid confusion, denote by $q_x$ the state
corresponding to $x \in X$. Note that this agreement does not
interfere with our earlier notation for sections, since $q_x =
\tau(q,x)=x$ in $A_{1+t}$.  For $p=\sum_{i=0}^\infty a_i t^i \in
R[[t]]$, we have (see Example~\ref{f})
\begin{gather*}
  q_x(p)= q_x \left( \sum_{i=0}^\infty a_it^i \right) =
  x + a_0 + (a_0+a_1)t + (a_1+a_2)t^2 + \dots = \\
  = x + \sum_{i=0}^\infty a_it^i + \sum_{i=0}^\infty a_it^{i+1} =
    x+(1+t)\sum_{i=0}^\infty a_it^i = x+(1+t)p.
\end{gather*}
Thus $q_x=\alpha^x\mu$, for $x \in X$, $\mu=q_0$,
$\alpha=q_1q_0^{-1}$, and therefore
\[ G(A_{1+t}) = \langle \ \{\ q_x \ |\  x \in X \ \}\rangle =
   \langle \alpha, \mu \rangle = G = L_n. \]
\end{proof}

\begin{proposition}
Let $g=a_0+a_1t+ \cdots + a_dt^d$ be a monic polynomial over
$R=Z/n\Z=X$ of degree $d \geq 1$, which is invertible in the power
series ring $R[[t]]$. The group of the automaton $A_g$ is the
lamplighter group
\[ L_{n,d} = (\oplus_\Z (\Z/n\Z)^d)\rtimes\Z  = (\Z/n\Z)^d \wr \Z, \]
\end{proposition}
\begin{proof}
This is just a straightforward generalization of the previous result.
First, note that, for $i= 0,\dots,d-1$ the maps $\alpha_i: R[[t]] \to
R[[t]]$ given by
\[ \alpha_i(p) = p + t^i \]
are tree automorphisms that have order $n$, commute, and generate a
copy of $(Z/n\Z)^d$. Let $G=\langle \alpha_0, \dots, \alpha_{d-1}, \mu
\rangle$, where $\mu:R[[t]] \to R[[t]]$ is the tree automorphism
given by
\[ \mu(p) = gp. \]
For $k \in \Z$ and $i=0,\dots,d-1$,
\[ \mu^k\alpha_i\mu^{-k}(p) = p + g^kt^i. \]
All these automorphisms have order $n$, commute, and generate the
normal closure $N$ of $\langle \alpha_0, \dots, \alpha_{d-1} \rangle$
in $G$, isomorphic to $\oplus_Z (\Z/n\Z)^d$. Moreover, since $\mu$ has
infinite order, we have $N \cap \langle \mu \rangle = 1$ and $G =
(\oplus_Z (\Z/n\Z)^d) \rtimes \Z \cong L_{n,d}$.

Let ${\bf y} = (y_0,\dots,y_{d-1}) \in X^d$ be a state of $A_g$. For
$p=\sum_{i=0}^\infty a_i t^i \in R[[t]]$, a straightforward
calculation shows that
\[ q_{\bf y}(p)= h_{\bf y} + gp, \]
where $h_{\bf y} = c_0 + c_1t + \cdots + c_{d-1}t^{d-1}$ and
\[ \begin{bmatrix} c_0 \\ c_1 \\ \vdots \\ c_{d-1} \end{bmatrix} =
   \begin{bmatrix} a_d & a_{d-1} & \dots & a_1 \\
                   0   & a_d     & \dots & a_2 \\
                   \vdots & \vdots & \ddots & \vdots \\
                   0   & 0       & \dots & a_d
   \end{bmatrix}
   \begin{bmatrix} y_0 \\ y_1 \\ \vdots \\ y_{d-1} \end{bmatrix}. \]
The above upper-triangular matrix is invertible over $R$ since its
determinant is 1 (recall that $a_d=1$). Therefore, for every
polynomial $h$ of degree smaller than $d$, there exists ${\bf
y}$ such that $q_{\bf y}(p) = h + gp$. In particular, $q_{\bf 0} =
\mu$ and
\[ G(A_g) = \langle \ \{\ q_{\bf y} \ |\ {\bf y} \in X^d \ \}\rangle =
   \langle \alpha_0, \alpha_1, \dots, \alpha_{d-1}, \mu \rangle = G = L_{n,d}. \]
\end{proof}

We now turn our attention to constructions of Baumslag-Solitar
groups and, more generally, ascending HNN extensions of free
abelian groups of finite rank.

\begin{proposition}\label{bs}
Let $m$ and $n$ be relatively prime integers greater than 1. The
group of the automaton $S_{m,n}$ is the Baumslag-Solitar
solvable group
\[ BS(1,m) = \langle \ a,t \ |\ tat^{-1} = a^m \ \rangle. \]
\end{proposition}

\begin{proof}
The infinite sequences over $X=\{0,1,\dots,n-1\}$ can be
interpreted as the elements of the ring $\Z_n$ of $n$-adic
integers. The state $s_0$ of the automaton $S_{m,n}$ simulates
the multiplication by $m$ in $\Z_n$. More generally the state
$s_i$, $i=0,\dots,m-1$ simulates the function $u \mapsto mu+i$
in $\Z_n$. Therefore the tree automorphism $\alpha=s_1s_0^{-1}$
is just the adding machine $u \mapsto u+1$ in $\Z_n$. Denote
$\mu=s_0$. Since $s_i = \alpha^i\mu$, we have $G(S_{m,n}) =
\langle \alpha, \mu \rangle$. Further,
\[ \mu\alpha\mu^{-1}(u) = \mu\alpha\left(\frac{1}{m}u\right) =
   \mu\left(\frac{1}{m}u+1\right) = u+m = a^m(u), \]
which shows that $\mu\alpha\mu^{-1} = \alpha^m$ in $G(S_{m,n})$
and therefore $G(S_{m,n})$ is a homomorphic image of $BS(1,m)$.
Both $\mu$ (multiplication by $m$ in $\Z_n$) and $\alpha$
(addition of 1 in $\Z_n$) have infinite order in $G(S_{m,n})$.
On the other hand, at least one of the images of $a$ and $b$
must have finite order in any proper homomorphic image of
$BS(1,m)$. Therefore $G(S_{m,n})=BS(1,m)$.
\end{proof}

The Baumslag-Solitar groups $BS(m,n)$ are subdivided as follows: if
$m=\pm1$ or $n=\pm1$, the group is solvable, and is realized by
automata, as explained above for $m=1$ or $n=1$. If $m=\pm n$, the
group is virtually $F_{|n|}\times\Z$, and therefore is realized by
automata, following e.g.~\cite{brunner-s:glnz}. Finally, if
$1\neq|m|\neq|n|\neq1$, then the group is not residually finite, so in
particular does not embed in the automorphism group of the rooted
tree, and thus cannot be realized by automata.

\begin{proposition}\label{G_M}
Let $M$ be an integer matrix of size $d \times d$ whose order is
infinite and the determinant
$m$ of $M$ is relatively prime to $n \geq 2$. There exists a
finite automaton on $n$ letters that defines the ascending $HNN$
extension
\[ G_M = \ \langle \ a_1, a_2, \dots , a_d, t \ || \ a_i \text{ commute },\
  a_i^t = a_1^{m_{1,i}}a_2^{m_{2,i}}\dots a_d^{m_{d,i}},\
  i=1,\dots,d \ \rangle ,\]
where $m_{i,j}$ is the entry in the row $i$ and column $j$ in
the matrix $M$.
\end{proposition}
\begin{proof}
Let $X$ be the alphabet $Y^d$, where $Y=\{0,1,\dots,n-1\}$ and
the elements of $X=Y^d$ are thought of as vector columns. The
infinite sequences over $X$ can be interpreted as the elements
of the free $\Z_n$-module of rank $d$, whose elements are also
considered as vector columns. Indeed, the free $\Z_n$-module of
rank $d$ consist of vector columns of size $d$ and each entry is
a member of $Z_n$, i.e., an infinite sequence over $Y$. Thus the
elements of the free module $\Z_n^d$ can be thought of as either
$d$-tuples of infinite sequences over $Y$ or as infinite
sequences of $d$-tuples over $Y$, i.e., infinite sequences over
$X$. The matrix $M$ is invertible over the ring $\Z_n$ since its
determinant $m$ is relatively prime to $n$. Thus we may think of
$M$ as being in $\GL_d(\Z_n)$, i.e., $M$ is a matrix of an
automorphism $\mu$ of the free module $\Z_n^d$ with respect to
the standard basis $({\bf e}_1,\dots,{\bf e}_d)$. Consider also,
for $i=1,\dots,d$, the translations $\alpha_i$ defined on
$\Z_n^d$ by  ${\bf u} \mapsto {\bf u} + {\bf e}_i$. Clearly, the
group generated by $\{\alpha_1,\dots,\alpha_d\}$ is the free
abelian group $\Z^d$. Moreover, for $i=1,\dots,d$,
\begin{gather*}
 \mu \alpha_i \mu^{-1}({\bf u}) = \mu \alpha_i (M^{-1}{\bf u}) =
 \mu(M^{-1}{\bf u} + {\bf e}_i) = \\
 = {\bf u} + M{\bf e}_i = {\bf u} + (m_{1,i},\dots,m_{d,i})^T =
 \alpha_1^{m_{1,i}} \cdots \alpha_d^{m_{d,i}}({\bf u}).
\end{gather*}
Thus the group $G=\langle \alpha_1,\dots,\alpha_d,\mu \rangle$
is a homomorphic image of the HNN extension $G_M$, under the
homomorphism that extends the map $t \mapsto \mu$, $a_i \mapsto
\alpha_i$, $i=1,\dots,d$. Under this homomorphism the image
$\langle \alpha_1,\dots,\alpha_d \rangle$ of $\langle a_1,\dots,a_d \rangle$
is free abelian group of rank d, the image $\langle \mu \rangle$ of
$\langle t \rangle$ is infinite cyclic group, and these two images
intersect trivially. However, in every proper homomorphic
image of $G_M$ the image of $\langle a_1,\dots,a_d \rangle$ is
not free abelian of rank $d$ or the image of $t$ has finite
order or these images have nontrivial intersection. This simply follows from
the fact that any non-trivial relation that can be added in $G_M$ must
have the form
\[ t^{k_0} = a_1^{k_1} \dots a_d^{k_d}, \]
where at least one of the integers $k_0,k_1,\dots,k_d$ is non-zero.
Thus the group $G$ is isomorphic to $G_M$.

The elements of $\Z_n^d$, being infinite sequences over $X$, can
be thought of as the boundary of the regular $n^d$-ary tree
$\tree$. It remains to be shown that there exists a finite
automaton, operating on $X$,  that defines a group of tree
automorphisms isomorphic to $G_M$. An example of such an
automaton is the automaton $T_{M,n}$ defined below, which
simulates the multiplication by the matrix $M$ in $\Z_n^d$.

More precisely, let
\[ \|M\| = \|M\|_\infty = \max_i \sum_j |m_{i,j}| \]
be the maximum absolute row sum norm (the max norm) of $M$
induced by the vector norm defined on vector columns ${\bf x}=
(x_1,\dots,x_d)^T$ by
\[ \|{\bf x}\|_\infty = \max_i |x_i|. \]
Let
\[ V = \{\ {\bf v} \ | \ {\bf v}=(v_1,\dots,v_d)^T \in \Z^d, \
  -\|M\| \leq v_i \leq \|M\| -1, \ i=1,\dots,d \ \}. \]
Define an automaton
\[ T_{M,n} = (T,X,\rho,\tau), \]
where $T= \{\ t_{\bf v} \ |\ {\bf v} \in V \ \}$ and $\rho : T
\times X \to X$ and $\tau: T \times X \to T$ are given by
\[ \rho(t_{\bf v},{\bf x})= (M{\bf x}+{\bf v}) \boxdot n \qquad\text{and}\qquad
   \tau(t_{\bf v},{\bf x})= t_{(M{\bf x}+{\bf v}) \div n}, \]
respectively, where $M{\bf x}+{\bf v}$ is calculated in $\Z^d$
and the remainder and quotient are defined by components.

The set of states is obviously finite (there are exactly
$(2\|M\|)^d$ states). Further, for ${\bf x} \in X$ and ${\bf v}
\in V$, the value of the $i$-th component of $M{\bf x} + {\bf
v}$ is between
\[ -\|M\|(n-1) - \|M\| = -\|M\| n
 \qquad\text{and}\qquad
    \|M\|(n-1) + \|M\| - 1 = \|M\| n -1, \]
respectively. This means that the $i$-th component in the
quotient $(M{\bf x}+{\bf v})\div n$ is always between $-\|M\|$
and   $\|M\| -1$ and therefore $t_{(M{\bf x}+{\bf v})\div n}$ is
always in $T$ and $\tau$ is well defined.

For fixed ${\bf v}$, the transformation ${\bf x} \mapsto (M{\bf x} +
{\bf v}) \boxdot n$ is a permutation of $X$ since the determinant $m$
of $M$ is relatively prime to $n$ (think of $X$ as the free module or
rank $d$ over the finite ring $\Z/n\Z$). Thus the automaton $T_{M,n}$
is invertible and each state defines an automorphism of the $n^d$-ary
tree $\Z_n^d$.

The state $t_{\bf v}$ defines the tree automorphism ${\bf u}
\mapsto M{\bf u}+{\bf v}$. Since $t_{{\bf e}_i} t_{{\bf
0}}^{-1}({\bf u}) = {\bf u} + {\bf e}_i$, the map $\alpha_i$ is
in $G(T_{M,n})$, for $i=1,\dots,d$. Finally, since $t_{{\bf
0}}=\mu$ we have
\[ G(T_{M,n}) = \langle \alpha_1,\dots,\alpha_d,t_{{\bf 0}} \rangle =
   G = G_M. \]
\end{proof}

Since every automaton group is a residually finite group with a word
problem that is solvable in exponential time, this shows that $G_M$
is always such a group. Note that the Dehn functions of the groups
$G_M$ have been carefully studied
(see for example~\cite{bridson-g:torus-bundles})
in the split case (i.e. when $M$ is in $\GL_n(\Z)$) and they are most
often exponential.

An analogous construction to the one above was used by Brunner and
Sidki in~\cite{brunner-s:glnz} to represent $\GL_n(\Z)$ by
automorphisms of the $2^n$-ary tree defined by finite automata.

\begin{example}
The automaton $T_{M,n}$ provided in the proof of
Proposition~\ref{G_M} is often not minimal automaton that
defines $G_M$. There is always a considerably smaller set of
states of $T_{M,n}$, closed under $\tau$, that defines a smaller
automaton and quite often still defines the same group. This
smaller automaton is defined as follows. Let $N_i$ and $P_i$ be
the sum of the negative entries and the positive entries,
respectively, in the row $i$ in $M$. Let
\[ V_S = \{ {\bf v} \ |\ {\bf v}=(v_1,\dots,v_d)^T \in \Z^d \ |\
        N_i \leq v_i \leq P_i-1,\ i=1,\dots,d \ \} \subset V. \]
Define an automaton
\[ S_{M,n} = (S,X,\rho,\tau), \]
where $S= \{\ s_{\bf v} \ |\ {\bf v} \in V_S \ \}$ and $\rho : S
\times X \to X$ and $\tau: S \times X \to S$ are defined as
restrictions of the maps in $T_{M,n}$. The minimal and the
maximal values of the $i$-th coordinate of $M{\bf x}+{\bf v}$,
for ${\bf x} \in X$ and ${\bf v} \in V_S$, are
\[ N_i(n-1)+N_i = N_i n \quad\text{and}\quad P_i(n-1)+P_i-1= P_in-1, \]
respectively, which means that $s_{(M{\bf x}+{\bf v}) \div n}$
is always in $S$ and the restriction $\tau$ is well defined.

For relatively prime $m,n \geq 2$ and $M=[m]$, the smaller automaton
$S_{M,n}$ is actually the automaton $S_{m,n}$ already defined before.
However, in general, the automaton $S_{M,n}$ does not generate $G_M$.
For example, if $M =
\begin{bmatrix} 3 & -1
\\ 0 & -1
\end{bmatrix}$, the automaton $S_{M,n}$ generates $BS(1,3)$,
while the larger automaton $T_{M,n}$ defines $G_M$. A sufficient
condition for the smaller automaton $S_{M,n}$ to generate $G_M$
is that the absolute row sum in each row of $M$ be at least $2$.
This condition is not necessary, as $S_{M,n}$, which has only 3
states, generates $G_M$ for $M=\begin{bmatrix} -2 & -1 \\ 1 & 0
\end{bmatrix}$.

Even when the smaller automaton $S_{M,n}$ generates $G_M$, there
sometimes exists yet smaller automata, operating on the same
alphabet, that define  $G_M$. For example, if $M=\begin{bmatrix}
-2 & 1 \\ 1 & -1 \end{bmatrix}$, the automaton $S_{M,n}$ has a
set of 6 states
\[ S = \{ \  s_{\bf v} \ |\ {\bf v} \in
  \{ \begin{bmatrix} 0 \\ 0 \end{bmatrix},
     \begin{bmatrix} -1 \\ 0 \end{bmatrix},
     \begin{bmatrix} -2 \\ 0 \end{bmatrix},
     \begin{bmatrix} 0 \\ -1 \end{bmatrix},
     \begin{bmatrix} -1 \\ -1 \end{bmatrix},
     \begin{bmatrix} -2 \\ -1 \end{bmatrix} \} \ \}. \]
However, the set of 4 states obtained from $S$ by exclusion of
the first and the last state is closed under the transition
function $\tau$ and is sufficient to generate $G_M$.
\end{example}


\section*{Dual automata}
The following general construction was considered before
in~\cite{macedonska-n-s:comm}:
\begin{definition}
Given an automaton $A=\A$, define the \emph{dual automaton} of
$A$, denoted by $A'$, by
\[ A' = (X,Q,\rho',\tau'), \]
where $\rho': X \times Q \to Q$ and $\tau' : X \times Q \to X$
are given by
\[ \rho'(x,q) = \tau(q,x) \qquad\text{and}\qquad
   \tau'(x,q) = \rho(q,x), \]
respectively.
\end{definition}

The definition of dual automaton confuses the letters with the states
and vice versa. In the graphical representation, for each edge
\[ q \overset{x|q(x)} \longrightarrow q_x \]
in the automaton $A$, there exist an edge
\[ x \overset{q|q_x} \longrightarrow q(x) \]
in the dual automaton $A'$. The confusion between states and letters
is possible because of the high symmetry present in the definition of
a finite transducer. In a sense, what we do is claim that not only
the states act on sequences of letters, but simultaneously the
letters act on sequences of states. We may say that, when the
automaton is in letter $x$ and reads the state $q$ it produces the
output state $q_x$ and lets the letter $q(x)$ handle the rest of the
input sequence of states. In other words, the domains of the
rewriting and transition functions are extended to arbitrary
sequences of states by
\begin{gather*}
 \rho(Wq,x) =  \rho(W,\rho(q,x))\\
 \tau(Wq,x) =  \tau(W,\rho(q,x))\tau(q,x),
\end{gather*}
for $x$ a letter in $X$, $q$ a state in $Q$ and $W$ a sequence
of states in $Q$. In the shorter notation, these equalities read
\[ Wq(x) = W(q(x)) \qquad\text{and}\qquad (Wq)_x = W_{q(x)}q_x. \]

\begin{definition}
An automaton $A$ is \emph{bi-invertible} if both $A$ and its
dual are invertible.
\end{definition}

It is easy to see that an automaton $A=\A$ is bi-invertible if,
for every state $q$ in $Q$, the restriction $\rho_q:X \to X$ is a
permutation of $X$, and, for every letter $x$ in $X$, the
restriction $\tau_x:Q \to Q$, given by $\tau_x(q)=\tau(q,x)$, is
a permutation of $Q$. The latter condition actually says that
the transition monoid (the transformation monoid over $Q$
generated by the maps $\tau_x:Q \to Q$, for $x \in X$) is a
group.

\begin{example}
For $n=2$, the automaton $\overline{A_{1+t}}$ from Example~\ref{barA}
is self-dual, i.e., it is isomorphic to its dual. Thus, somewhat
trivially, $\overline{A_{1+t}}$ is bi-invertible. The transition
monoid is the cyclic group of order $2$.

For $n \geq 3$ the  automaton $\overline{A_{1+t}}$ is also
bi-invertible. Identify the set of states $\overline{X}$ with $X$.
The letter $x$ in $X$ induces the permutation $y \mapsto -y+x$ on the
set of states. The transition monoid is then the subgroup of
permutations of the state set $X$ generated by the permutations $y
\mapsto -y+x$, for $x \in X$. This group is generated by the two
involutions $y \mapsto -y$ and $y \mapsto -y+1$. Since their
composition is a cyclic permutation of order $n$, the transition
monoid of the bi-invertible automaton $\overline{A_{1+t}}$ is the
dihedral group $D_n$ (the symmetry group of the regular $n$-gon).
\end{example}

\begin{proposition}
Let $m,n \geq 2$ be relatively prime integers, $m'$ an integer that is a
multiplicative inverse of $m$ modulo $n$ and $n'$ an integer that is a
multiplicative inverse of $n$ modulo $m$. Define the automaton
\[ D_{m,n} = (D,X,\rho,\tau), \]
where $D=\{d_0,\dots,d_{m-1}\}$, $X = \{x_0,\dots,x_{n-1}\}$ and
$\rho: D \times X \to X$ and $\tau : D \times X \to D$ are given
by
\[ \rho(d_i,x_j)= x_{m'(j-i) \boxdot n} \qquad\text{and}\qquad
   \tau(d_i,x_j)= d_{n'(i-j) \boxdot m}, \]
respectively.
\begin{itemize}
\item[(a)] The definition of the automaton $D_{m,n}$ does not
depend on the choice of $m'$ and $n'$.
\item[(b)] The automaton $D_{m,n}$ is the inverse of the automaton $S_{m,n}$.
\item[(c)] The dual of the automaton $D_{m,n}$ is $D_{n,m}$.
\item[(d)] The automaton $D_{m,n}$ is bi-invertible.
\item[(e)] The group $G(D_{m,n})$ is the Baumslag-Solitar group
$BS(1,m)$.
\end{itemize}
\end{proposition}
\begin{proof}
(a) Clear.

(b) Consider the quantities
\[ y_{i,j} = n ( n'(i-j) \boxdot m)+j \qquad\text{and}\qquad
   z_{i,j}= m ( m'(j-i) \boxdot n) +i, \]
for $i=0,\dots,m-1$ and $j=0,\dots,n-1$. Since
\[ y_{i,j} \boxdot n = j = z_{i,j} \boxdot n \quad\text{and}\qquad
   y_{i,j} \boxdot m = i = z_{i,j} \boxdot m, \]
the quantities $y_{i,j}$ and $z_{i,j}$ differ by a multiple of
$mn$, according to the Chinese Remainder Theorem. However, $0
\leq y_{i,j},\ z_{i,j} \leq mn-1$ and therefore $y_{i,j}
=z_{i,j}$. Thus
\[ (m (m'(j-i) \boxdot n)+i) \div n = z_{i,j} \div n = y_{i,j} \div n =
    n'(i-j) \boxdot m, \]
which shows that the automaton $D_{m,n}$ is just the automaton
$\overline{S_{m,n}}$ in disguise.

(c) Evident from the symmetry in the definition of $D_{m,n}$.

(d) It follows from (b) that $D_{m,n}$ is invertible, and then from (c) that
it is bi-invertible.

(e) Every invertible automaton generates the same group as its
inverse automaton, so the result follows from (b) and
Proposition~\ref{bs}.
\end{proof}

\begin{figure}[!hbt]
\begin{center}
\includegraphics{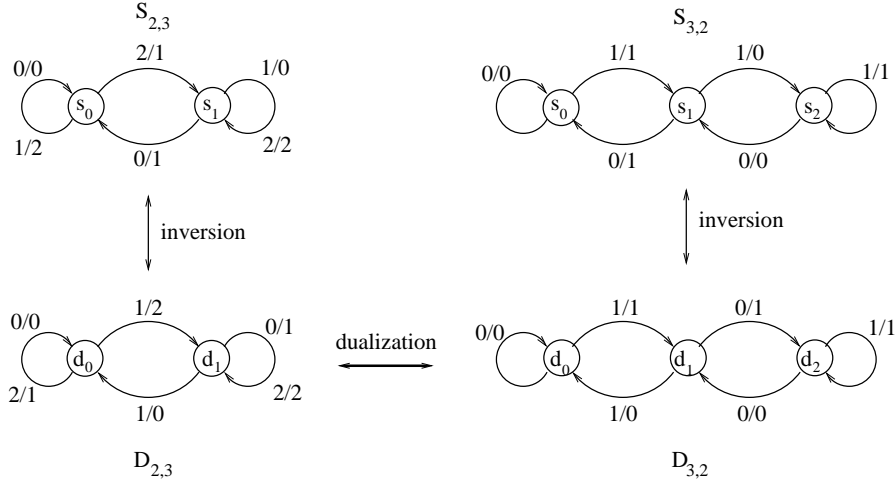}
\end{center}
\caption{Relations between the automata $S_{2,3}$, $S_{3,2}$,
$D_{2,3}$ and $D_{3,2}$ \label{sd23}}
\end{figure}

The above proposition says that the automata $S_{m,n}$, $S_{n,m}$, $D_{m,n}$ and
$D_{n,m}$ are related as follows
\[ S_{m,n} \overset{inversion}{\longleftrightarrow} D_{m,n}
   \overset{dualization}{\longleftrightarrow}
   D_{n,m} \overset{inversion}{\longleftrightarrow} S_{n,m}. \]
These relations are depicted in Figure~\ref{sd23} for $m=3$ and $n=2$.

The above relations show that there is an interesting connection
between $BS(1,m)$ and $BS(1,n)$ for any pair of relatively prime
integers greater $m,n \geq 2$. Indeed, $BS(1,m)$ is defined by the
automaton $D(m,n)$ having $m$ states and operating on an $n$-letter
alphabet, while $BS(1,n)$ is defined by the automaton $D_{n,m}$ on
$n$ states operating on an $m$-letter alphabet, and the latter
automaton is obtained by simple dualization procedure that
``confuses'' states with letters and the other way around in the
former automaton.


\section*{Composition of automata and Schreier graphs}
\begin{proposition}\label{minus}
The automaton $T_{M,n}$ can be obtained from the automaton $T_{-M,n}$
(and vice versa) by multiplying on the left each permutation
$\rho_q:X \to X$, for $q$ a state in $T_{-M,n}$, by the involution
${\bf x} \mapsto (-{\bf x} - {\bf 1}) \boxdot n$, where ${\bf 1} =
\sum_{i=1}^d{\bf e_i}$. In exactly the same way $S_{M,n}$ can be
obtained from $S_{-M,n}$.
\end{proposition}
\begin{proof}
Note that $\|M\|=\|-M\|$. Thus the set of vectors $V$ used to
index the states in $T_{M,n}$ and $T_{-M,n}$ is the same. In
order to avoid confusion, change the names of the states of
$T_{-M,n}$ to $k_{\bf v}$, ${\bf v} \in V$. Consider the
bijection $f$ between the states of $T_{-M,n}$ and the states of
$T_{M,n}$ given by $k_{\bf v} \mapsto t_{-{\bf v} - {\bf 1}}$.
Then,
\[ f(\tau(k_{\bf v},{\bf x})) = f(k_{(-M{\bf x}+{\bf v})\div n}) =
   t_{ - (-M{\bf x}+{\bf v})\div n - {\bf 1} } \]
and
\[ \tau(f(k_{\bf v}),{\bf x})= \tau(t_{-{\bf v} - {\bf 1}}, {\bf x}) =
    t_{(M{\bf x} -{\bf v} - {\bf 1})\div n}, \]
for ${\bf v} \in V$ and ${\bf x} \in X$. One can easily verify
that $(a-1) \div n = -(-a \div n) -1$ for any integer $a$. Thus
$f(\tau(k_{\bf v},{\bf x})) = \tau(f(k_{\bf v}),{\bf x})$, which
means that $f$ is compatible with the transition functions
defined in the two automata, i.e., the transition in the
automaton $T_{-M,n}$ at $k_{\bf v}$ behaves exactly as the
transition in $T_{M,n}$ at $f(k_{\bf v})$.

Let $\xi:X \to X$ be the involution ${\bf x} \mapsto (-{\bf x} - {\bf
1}) \boxdot n$. Then, for ${\bf v} \in V$ and ${\bf x} \in X$,
\[ \xi(\rho(k_{\bf v},{\bf x})) = \xi( (-M{\bf x}+{\bf v}) \boxdot n) =
   (M{\bf x} - {\bf v} - {\bf 1}) \boxdot n =
   \rho(t_{-{\bf v}-{\bf 1}},{\bf x}) = \rho(f(k_{\bf v}),{\bf x}). \]
This proves the first claim. Note that if $\rho(k_{\bf v},{\bf
x})$ were equal to $\rho(f(k_{\bf v}),{\bf x})$, then  $f$ would
have been an isomorphism between the two automata.

The second claim follows easily, since $f$ maps bijectively the
states of $S_{-M,n}$ onto the states of $S_{M,n}$.
\end{proof}

The way in which $T_{M,n}$ is obtained from $T_{-M,n}$ is just a
special case of a more general construction of composition of
automata. Informally, given two automata $A$ and $B$ operating over
the same alphabet $X$ one wants to construct an automaton that
operates over the same alphabet and, for every pair of states $p$ and
$q$ in $A$ and $B$, respectively, contains a state that acts on a
word $w$ over $X$ exactly as $p$ would act on the output of the
action of $q$ on $w$ (i.e., it acts as $q$ followed by $p$).

\begin{definition}
Let $A=(P,X,\rho_2,\tau_2)$ and $B=(Q,X,\rho_1,\tau_1)$ be two finite
automata. The \emph{composition} of the two automata, denoted $AB$,
is the automaton
\[ AB = (P \times Q, X, \rho, \tau) \]
where $\rho: (P \times Q) \times X \to X$ and $\tau: (P \times
Q) \times X \to P \times Q$ are given by
\[ \rho((p,q),x) = \rho_2(p,\rho_1(q,x)) \qquad\text{and}\qquad
   \tau((p,q),x) = (\tau_2(p,\rho_1(q,x)),\tau_1(q,x)) , \]
respectively.
\end{definition}

It is easy to verify that the composition of two invertible
automata as above is an invertible automaton in which
\[ \rho_{(p,q)} = \rho_p \rho_q \qquad\text{and}\qquad
    (p,q)_x = (p_{q(x)},q_x), \]
for $p$ a state in $P$, $q$ a state in $Q$ and $x$ a letter in $X$.
The above equalities are consistent with (\ref{composition}),
indicating that the state $(p,q)$ in $AB$ defines the composition
$pq$ of the tree automorphisms $p$ and $q$.

\begin{example}
Consider again, as in Proposition~\ref{minus}, the relation between
$T_{-M,n}$ and $T_{M,n}$. The automaton $A$ on a single state $q$,
for which $\rho_q$ is the permutation $\xi:{\bf x} \mapsto (-{\bf x}
- {\bf 1}) \boxdot n$, defines the cyclic group of order $2$. The
automorphism $q$ of the $n^d$-ary tree defined by $q$ is the
involution ${\bf u} \mapsto -{\bf u}-{\bf 1}$. The composition
$AT_{-M,n}$ is isomorphic to $T_{M,n}$ under the correspondence
$(q,k_{\bf v}) \leftrightarrow t_{-{\bf v} - {\bf 1}}$.
\end{example}

In the light of the observation that the state $(p,q)$ is the
composition of invertible automata $A$ and $B$ represents the
composition of the tree automorphisms represented by $p$ and $q$, the
following remark is obvious.

\begin{proposition}
Let $A=\A$ be an invertible automaton. The group $G(A^k)$ of the
automaton $A^k$ is the subgroup of $G(A)$ generated by all words of
length $k$ over the states of $A$.
\end{proposition}

\begin{proposition}\label{sm1m2}
Let $m$, $m_1$, $m_2$ and $n$ be positive integers such that $m$,
$m_1$ and $m_2$ are all relatively prime to $n$, and let $k \geq 1$.
Then $ G(S_{m_1,n} S_{m_2,n}) = G(S_{m_1m_2,n}) = BS(1,m_1m_2)$ and
$G((S_{m,n})^k) = G(S_{m^k,n}) = BS(1,m^k)$ Moreover,
\[ S_{m_1,n} S_{m_2,n} = S_{m_1m_2,n} \qquad\text{and}\qquad
   S_{m,n}^k = S_{m^k,n}. \]
\end{proposition}
\begin{proof}
All claims follow from the fact that $S_{m_1,n} S_{m_2,n} \cong
S_{m_1m_2,n}$. The latter can be easily proved by observing that an
automaton isomorphism (fixing the alphabet) from $S_{m_1,n}
S_{m_2,n}$ to $S_{m_1m_2,n}$ is given by
\[ (s_i, s_j) \mapsto s_{m_1j+i}, \]
for $i \in \{0,\dots,m_1-1\}$, $j \in \{0,\dots,m_2-1\}$.
\end{proof}

\begin{proposition}\label{overlineAB}
For any two invertible automata $A=(P,X,\rho_2,\tau_2)$ and $B=(Q,X,\rho_1,\tau_1)$,
the automaton $AB$ is invertible and
\[ \overline{AB} = \overline{B} \ \overline{A}. \]
More generally, for any invertible automata $A_1,\dots,A_k$ over the
same alphabet, the automaton $A_1 \dots A_k$ is invertible and
\[ \overline{A_1 \dots A_k} = \overline{A_k}  \dots \overline{A_1}. \]
\end{proposition}
\begin{proof}
The automaton $AB$ is invertible since, for $(p,q)$ a state in $AB$,
the map $\rho_{(p,q)}:X \to X$ is invertible. The latter is clear
since $\rho_{(p,q)}$ is the composition $\rho_p\rho_q$ of invertible
maps.

Consider the edge
\[(p,q) \overset{x|p(q(x))}{\longrightarrow} (p_{q(x)},q_x) \]
in $AB$ and its corresponding edge
\begin{equation} \label{AB}
\overline{(p,q)} \overset{p(q(x))|x}{\longrightarrow}
\overline{(p_{q(x)},q_x)}
\end{equation}
in $\overline{AB}$. Let $y = p(q(x))$ and consider the edge
\begin{equation} \label{BA}
   (\overline{q},\overline{p})
   \overset{y|\overline{q}(\overline{p}(y))}{\longrightarrow}
   (\overline{q}_{\overline{p}(y)},\overline{p}_y)
\end{equation}
in $\overline{B} \ \overline{A}$. We have
\begin{gather*}
\overline{q}(\overline{p}(y)) = q^{-1}(p^{-1}(y)) = x, \\
\overline{q}_{\overline{p}(y)} = \overline{q}_{p^{-1}(y)} = \overline{q}_{q(x)} =
\overline{q_{q^{-1}(q(x))}} = \overline{q_x}, \\
\overline{p}_y = \overline{p_{p^{-1}(y)}} = \overline{p_{q(x)}}.
\end{gather*}
Thus the edge~(\ref{BA}) can be rewritten as
\begin{equation}\label{BA'}
  (\overline{q},\overline{p}) \overset{p(q(x))|x}{\longrightarrow}
   (\overline{q_x},\overline{p_{q(x)}}).
\end{equation}
The canonical bijection $\overline{(p,q)} \mapsto
(\overline{q},\overline{p})$ maps the edge~(\ref{AB}) to the
edge~(\ref{BA'}). Thus $\overline{AB}$ and $\overline{B} \
\overline{A}$ are canonically isomorphic.
\end{proof}

\begin{proposition}\label{dm1m2}
Let $m$, $m_1$, $m_2$ and $n$ be positive integers such that $m$,
$m_1$ and $m_2$ are all relatively prime to $n$, and let $k \geq 1$.
Then $ G(D_{m_2,n} D_{m_1,n}) = G(D_{m_1m_2,n}) = BS(1,m_1m_2)$ and
$G((D_{m,n})^k) = G(D_{m^k,n}) = BS(1,m^k)$ Moreover,
\[ D_{m_2,n} D_{m_1,n} = D_{m_1m_2,n} \qquad\text{and}\qquad
   D_{m,n}^k = D_{m^k,n}. \]
\end{proposition}
\begin{proof}
This is a direct corollary of Proposition~\ref{sm1m2} and
Proposition~\ref{overlineAB}. The only point worth mentioning is that
the canonical isomorphism from $D_{m_2,n} D_{m_1,n}$ to
$D_{m_1m_2,n}$, which is composed from the two canonical isomorphisms
in Proposition~\ref{sm1m2} and Proposition~\ref{overlineAB} is given
by
\[ (d_j, d_i) \mapsto d_{m_1j+i}, \]
for $i \in \{0,\dots,m_1-1\}$, $j \in \{0,\dots,m_2-1\}$. Indeed,
\[ (d_j,d_i) = (\overline{s_j},\overline{s_i}) \mapsto
\overline{(s_i,s_j)} \mapsto \overline{s_{m_1j+i}} = d_{m_1j+i}. \]
\end{proof}

Consider an invertible automaton $A=\A$. The action of the group
$G(A)$ on the $k$-th level of the tree $X^*$ can be depicted by a
finite graph, known as the \emph{Schreier graph} of the action, as
follows. The vertices are the $k$-letter words over $X$ and, for each
vertex $u=x_1x_2\dots x_k$ and a generator (state) $q$ in $Q$, a
directed edge labelled by $q$ connects $u$ to $q(u)$. In our
situation we can enrich the structure of this graph by labelling the
edge from $u$ to $q(u)$ by $q|q_u$. With this the Schreier graph
becomes the graphical representation of an automaton. Denote the
obtained automaton by $Sch_k(A)$ and call it the $k$-level Schreier
automaton of $A$. For $k=1$, the obtained Schreier automaton is just
the dual automaton $A'$, i.e.,
\[ Sch_1(A) = A' \]

\begin{proposition}\label{sch}
Let $\A$ be an invertible automaton. Then, for all positive
integers $k$,
\[ Sch_k(A) \cong (A')^k, \]
where the isomorphism canonically maps the $k$-letter word
$u=x_1\dots x_k$ over $X$ (a state in $Sch_k(A)$) to the state
$(x_k,\dots,x_1)$ in $(A')^k$.
\end{proposition}
\begin{proof}
It is clear that the canonical map is bijection between the states of
$Sch_k(A)$ and $(A')^k$.

Let $u=x_1 \dots x_k$ be an arbitrary word over $X$ and $q$ a state
in $A$. The edges
\begin{gather*}
  x_1 \overset{q|q_{x_1}}{\longrightarrow} q(x_1), \\
  x_2 \overset{q_{x_1}|q_{x_1x_2}}{\longrightarrow} q_{x_1}(x_2), \\
  \dots, \\
  x_k \overset{q_{x_1 \dots x_{k-1}}|q_{x_1 \dots x_k}}{\longrightarrow} q_{x_1\dots x_{k-1}}(x_k)
\end{gather*}
in $A'$ imply that the edge corresponding to $(x_k, \dots ,x_1)$ and
$q$ in $(A')^k$ is
\[ (x_k, \dots, x_1) \overset{q | q_u}{\longrightarrow}
   (q_{x_1 \dots x_{k-1}}(x_k), \dots, q_{x_1}(x_2),q(x_1)). \]
Since $q(x_1)q_{x_1}(x_2) \dots q_{x_1 \dots x_{k-1}}(x_k) =
q(x_1\dots x_k)=q(u)$, the corresponding edge in $Sch_k(A)$ is
\[ u \overset{q | q_u}{\longrightarrow} q(u), \]
so the result follows.
\end{proof}

Thus, in general, the $k$-fold power of the dual graph of $A$ looks
exactly the same as the Schreier graph of the action of $A$ on level
$k$, with the only difference being the reversal in the order in the
$k$-tuples representing the states of these two automata.

\begin{proposition}
For relatively prime integers $m,n \geq 2$ and $k \geq 1$,
\[ BS(1,n^k) = G(D_{n^k,m}) = G((D_{n,m})^k) = G(Sch_k(D_{m,n})).\]
Moreover
\[ D_{n^k,m} \cong (D_{n,m})^k \cong Sch_k(D_{m,n}).\]
\end{proposition}
\begin{proof}
First identify the symbol $d_i$ for the states in all automata above
with the symbol $i$.

By Proposition~\ref{dm1m2} the automaton $D_{n^k,m}$ looks exactly
the same as $(D_{n,m})^k$, except that the state $i$ in $D_{n^k,m}$
corresponds to the $k$-tuple $(a_{k-1},\dots,a_0)$, where
\[ i = a_{k-1} n^{k-1} + \cdots + a_1 n  + a_0 \]
is the $k$-digit $n$-ary representation of the non-negative integer
$i$, for $i=0,\dots,n^k-1$.

By Proposition~\ref{sch} the $k$-level Schreier automaton
$Sch_k(D_{m,n})$ looks also exactly the same as $(D_{n,m})^k$, with
the state $(a_{k-1},\dots,a_0)$ in $(D_{n,m})^k$ corresponding to
$a_0a_1\dots a_{k-1}$ in $Sch_k(D_{m,n})$.
\end{proof}

\begin{example}
Figure~\ref{d43} depicts the automaton $D_{4,3}$ and illustrates the
previous proposition.

\begin{figure}[!hbt]
\begin{center}
\includegraphics{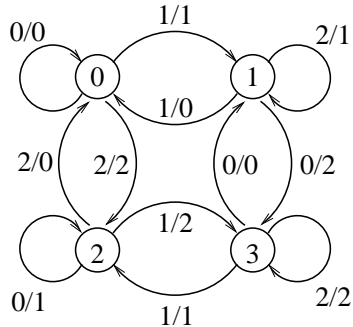}
\end{center}
\caption{The automaton $D_{4,3}$ \label{d43}}
\end{figure}

The square automaton $(D_{2,3})^2$ looks exactly the same as
$D_{4,3}$, except that the state 0 corresponds to the pair $(0,0)$,
the state 1 to the pair $(0,1)$, the state 2 to the pair $(1,0)$ and
the state 4 to the pair $(1,1)$

The second level Schreier automaton of $D_{3,2}$ also looks exactly
the same as $D_{4,3}$, except that 0 corresponds to 00, 1 to 10, 2 to
01 and 3 to 11.
\end{example}

\section*{Acknowledgments}
The second author would like to thank Nata{\v s}a Jonoska, Mile
Kraj{\v c}evski and the Department of Mathematics at University of
South Florida for their hospitality during my extended visit
during which most of the manuscript was completed.

Thanks to the referee for his/her help in improving the
presentation.

\bibliographystyle{alpha}

\def\cprime{$'$}

\end{document}